\documentclass[a4paper,12pt]{article}

\usepackage{hyperref}
\usepackage[pdftex]{graphicx}
\usepackage[section] {placeins}
\usepackage{enumerate}
\usepackage[normalem]{ulem}
\usepackage{amsmath}
\usepackage{amsthm}
\usepackage{latexsym}
\usepackage{amsfonts}
\usepackage{amssymb}
\usepackage{mathtools}
\usepackage{longtable}

\usepackage{tikz-cd}
\usetikzlibrary{graphs}
\usetikzlibrary{arrows.meta}
\usetikzlibrary{shapes.geometric}

\usepackage{authblk}

\usepackage{latexsym,amssymb,amsfonts, amsthm, amsmath}
\usepackage[english]{babel}
\usepackage{bm}
\usepackage{mathrsfs}
\usepackage{algpseudocode}
\usepackage{algorithm}

\algnewcommand{\IIf}[1]{\State\algorithmicif\ #1\ \algorithmicthen}
\algnewcommand{\EndIIf}{\unskip\ \algorithmicend\ \algorithmicif}

\newtheorem{thm}{Theorem}[section]
\newtheorem{lem}[thm]{Lemma}
\newtheorem{cor}[thm]{Corollary}
\newtheorem{prop}[thm]{Proposition}

\newtheorem{rem}[thm]{Remark}
\newtheorem*{thm*}{Theorem \ref{th:main3}}
\newtheorem*{thm*2}{Theorem \ref{th:main1}}

\theoremstyle{definition}
\newtheorem{defn}[thm]{Definition}

\newcommand{\Z}{\mathbb{Z}}

\def\Dic{\mathrm{Dic}}

\usepackage{fullpage}

\makeatletter
\newcommand{\subjclass}[2][2010]{%
\let\@oldtitle\@title%
\gdef\@title{\@oldtitle\footnotetext{#1 \emph{Mathematics subject classification.} #2.}}%
}
\newcommand{\keywords}[1]{%
\let\@@oldtitle\@title%
\gdef\@title{\@@oldtitle\footnotetext{\emph{Keywords and phrases.} #1.}}%
}
\makeatother

\begin{document}
\title{Weak Freiman isomorphisms \\and sequencings of small sets}
\author[1]{Simone Costa}
\author[2]{Stefano Della Fiore}

\affil[1]{DICATAM, Sez.~Matematica, Universit\`a degli Studi di Brescia, Via Branze~43, I~25123 Brescia, Italy}
\affil[2]{DI, Universit\`a degli Studi di Salerno, Via Papa Giovanni Paolo II~132, \newline I~84084 Fisciano, Italy, }
\maketitle
\begin{abstract}
In this paper, we introduce a weakening of the Freiman isomorphisms between subsets of non necessarily abelian groups.

Inspired by the breakthrough result of  Kravitz, \cite{NK}, on cyclic groups, as a first application, we prove that any subset of size $k$ of the dihedral group $D_{2m}$ (and, more in general, of a class of semidirect products) is sequenceable, provided that the prime factors of $m$ are larger than $k!$.  Also, a refined bound of $k!/2$ for the size of the prime factors of $m$ can be obtained for cyclic groups $\mathbb{Z}_m$, slightly improving the result of \cite{NK}.
Then, applying again the concept of weak Freiman isomorphism, we show that any subset of size $k$ of the dicyclic group $\Dic_{m}$ is sequenceable, provided that the prime factors of $m$ are larger than $k^k$.
\end{abstract}
\section{Introduction}
A subset of a group is {\em sequenceable} if there is an ordering~$(x_1, \ldots, x_k)$ of its elements such that the partial sums~$(s_0, s_1, \ldots, s_k)$, given by $s_0 = 0$ and $s_i = \sum_{j=1}^i x_i$ for $1 \leq i \leq k$, are distinct, with the possible exception that we may have~$s_k = s_0 = 0$; such orderings will be called {\em sequencings} of the subset. Note that, even if we are using the additive notation, we are not assuming that the group is abelian.

Starting from Graham and Erd\H{o}s (see \cite{GR} and \cite{EG}), several conjectures and questions concerning the sequenceability of subsets of abelian groups appear in the literature, which have been combined and summarized (see \cite{AL20} and \cite{CDOR}) into the conjecture that if a subset of an abelian group does not contain~0 then it is sequenceable. For related conjectures and applications of these results, we also refer to the papers \cite{AKP} and \cite{ADMS16}, and to the surveys \cite{PD} and \cite{OllisSurvey}.
Even though these papers are more focused on abelian groups, also the case of dihedral ones (and, more in general, semidirect products) have been studied by Ollis in \cite{OllisDi} and by Costa, Della Fiore, and Ollis in \cite{CDO}.

In \cite{NK}, Kravitz, prove the breakthrough result that any subset $S$ of size $k$ of $\mathbb{Z}_m\setminus\{0\}$ is sequenceable provided that the prime factors of $m$ are larger than $(k-1)^{k+1}$. Before his paper, the literature was concentrated on cases where either the group or the set is smaller than a given constant (see, for instance, \cite{HOS19}, \cite{CDOR} or \cite{CMPP18}) or where the partial sums were compared only locally (see \cite{CD} and \cite{CD1}). Kravitz, in his paper, applies the so-called rectification method: he embeds the elements of $S$ in $\mathbb{Z}$ through a Freiman isomorphism and then he proves its sequenceability result on $\mathbb{Z}$.

Here we show that by considering a weakening of the Freiman isomorphism method it is possible to generalize this procedure to obtain the same result for the dihedral group $D_{2m}$ (and, more in general, for a class of semidirect products) provided that the prime factors of $m$ are larger than $k!$.
Also,  a refined bound of $k!/2$ for the size of the prime factors of $m$ can be obtained for cyclic groups $\mathbb{Z}_m$, slightly improving the result of \cite{NK}.

Then, following the line of the proof of Theorem 3.1 of \cite{BLR} we show that any subset of size $k$ of the dicyclic group $\Dic_m$ is weak Freiman isomorphic to a subset of the semidirect product $\mathbb{Z}\rtimes_{\varphi} \mathbb{Z}_4$ (where $\varphi(1)=-id$) provided that the prime factors of $m$ are larger than $k^k$. It will follow that, under the same assumption, any subset of size $k$ of the dicyclic group $\Dic_{m}$ is sequenceable.
\section{Weak Freiman isomorphisms }
The rectification method used by Kravitz in \cite{NK} relies on the concept of Freiman isomorphism. Following \cite{TV1} here, and in all the paper, we will use the additive notation even if the groups are not assumed to be abelian (if not specified).
We recall the following definition (see Tao and Vu, \cite{TV1}).
\begin{defn}
Let $\ell\geq 1$, and let $S$, $S'$ be sets with ambient groups $(V,+_V)$ and $(W,+_W)$ respectively. A Freiman homomorphism of order
$\ell$, say $\phi$, from $(S, V)$ to $(S',W)$ (or more succinctly from $S$ to $S'$) is a map $\phi : S \rightarrow S'$ with the property that:
$$a_1+_V a_2+_V  \dots+_V a_{\ell}=a_1'+_V a_2'+_V\dots +_V a_{\ell}'$$ implies $$\phi(a_1) +_W\phi(a_2)+_W\dots+_W\phi(a_{\ell})=\phi(a_1') +_W\phi(a_2') +_W\dots+_W\phi(a_{\ell}').$$
If in addition there is an inverse map $\phi^{-1}: S' \rightarrow S$ which is also a Freiman homomorphism of order $\ell$ from $(S',W)$ to $(S, V)$, then we say that $\phi$ is a Freiman isomorphism of order $\ell$ and that $(S, V)$ and $(S',W)$ are Freiman isomorphic of order $\ell$.
\end{defn}
In case $V=\mathbb{Z}_m$ and $W=\mathbb{Z}$, it is knowns that (see \cite{Lev}) such isomorphism exists provided that the prime factors of $m$ are large enough (i.e. $\ell^{k}$).

Inspired by this concept, and to adapt the rectification procedure to the case of semidirect products, we introduce the following weakening.
\begin{defn}
Let $k\geq 1$, and let $S$, $S'$ be sets with ambient groups $(V,+_V)$ and $(W,+_W)$ respectively. A weak Freiman homomorphism of order
$k$, say $\phi$, from $(S, V)$ to $(S',W)$ (or more succinctly from $S$ to $S'$) is a map $\phi : S \rightarrow S'$ with the property that, for any ordered subset $(a_1,\dots,a_{\ell})$ of cardinality $\ell\leq k$ of $A$:
$$a_1+_Va_2+_V\dots +_V a_{\ell}=0\rightarrow \phi(a_1)+_W\phi(a_2)+_W \dots+_W\phi(a_{\ell})=0.$$
If in addition there is an inverse map $\phi^{-1}: S' \rightarrow S$ which is also a weak Freiman homomorphism of order $k$ from $(S',W)$ to $(S, V)$, then we say that $\phi$ is a weak Freiman isomorphism of order $k$ and that $(S, V)$ and $(S',W)$ are Freiman weakly isomorphic of order $k$.

In case $k=|S|$, we will simply say that $\phi$ is a weak Freiman homomorphism (resp. isomorphism).
\end{defn}
Roughly speaking a weak Freiman homomorphism only preserves the sets which sum to zero in a given order. A Freiman homomorphism instead preserves all the equality between sequences sums and in particular it allows to have repeated elements in the sum. Note that, if we assume $0\in S$, a Freiman homomorphism also preserves the sequences that sum to zero.

In the following, we consider groups of the form $(\mathbb{Z}_m\rtimes_{\varphi} H,+)$ where $H$ is finite, $(H,+_H)$ is not necessarily abelian and where $\varphi$ is a group homomorphism $\varphi: H\rightarrow Aut(\mathbb{Z}_m)$. We will also assume that the homomorphism $\varphi$ satisfies the following property:
\begin{itemize}
\item[$(*)$]\  given $h\in H, \varphi(h)$ is either the identity or $\varphi(h)(x)=-x$ for any $x\in \mathbb{Z}_m$.
\end{itemize}
Note that the latter is an automorphism of $\mathbb{Z}_m$.
\begin{rem}
Two classes of semidirect product that satisfy property $(*)$ are the dihedral groups and the standard direct product.
Indeed:
\begin{itemize}
\item[1)] The dihedral group $D_{2m}$ (called also $Dih(\mathbb{Z}_m)$) can be defined as $\mathbb{Z}_m\rtimes_{\varphi} \Z_2$ where
$$\varphi(0)=id \mbox{ and }\varphi(1)(x)=-x.$$
This group can be defined also on $\mathbb{Z}$. In this case, we will denote it by $Dih(\mathbb{Z})$.
\item[2)] Any direct product $\mathbb{Z}_m\times H$ can be defined as $\mathbb{Z}_m\rtimes_{\varphi} H$ where
$$\varphi(h)=id \mbox{ for any }h\in H.$$
\end{itemize}
\end{rem}

Here, following the proof of Theorem 22 of \cite{CDOR} (see also Theorem 3.14 of \cite{CDO}, and Theorem 3.1 of \cite{CP20} where an explicit bound was not provided), we prove the following theorem.
\begin{thm}\label{main}
Let $H$ be a finite group, let $\varphi$ satisfy property $(*)$. Then any subset $S$ of size $k$ of $\mathbb{Z}_m\rtimes_{\varphi} H\setminus\{0\}$ is Freiman weakly isomorphic to a subset $S'$ of size $k$ of  $\mathbb{Z}\rtimes_{\varphi} H\setminus\{0\}$ provided that the prime factors of $m$ are all greater than $k!$.
\end{thm}
First, we need to recall some elementary facts about the solutions of a linear system over fields and commutative rings.
\begin{thm}[Corollary of Cramer's Theorem]\label{Cramer}
Let
\begin{equation}\label{sys1}\begin{cases}
m_{1,1}x_1+\dots+ m_{1,k}x_k=b_1\\
\dots\\
m_{l,1}x_1+\dots+ m_{l,k}x_k=b_l
\end{cases} \end{equation}
be a linear system over a commutative ring $R$ whose associated matrix $M=(m_{i,j})$ is an $l\times k$ matrix over $R$. Let $(x_1,\dots, x_k)$ be a solution of the system and let $M'$ be a $k'\times k'$ square, nonsingular (i.e.~$\det(M')$ is invertible in $R$), submatrix of $M$. We assume, without loss of generality, that $M'$ is obtained by considering the first $k'$ rows and the first $k'$ columns of $M$. Then $x_1,\dots, x_k$ satisfy the following relation

$$\begin{bmatrix}
x_{1} \\
x_{2} \\
\vdots \\
x_{k'}
\end{bmatrix} = \begin{bmatrix}
m_{1,1} & m_{1,2} & m_{1,3} & \dots & m_{1,k'} \\
m_{2,1} & m_{2,2} & m_{2,3} & \dots & m_{2,k'} \\
\vdots & \vdots & \vdots & & \vdots \\
m_{k',1} & m_{k',2} & m_{k',3} & \dots & m_{k',k'}
\end{bmatrix} ^{-1} \begin{bmatrix}
b_1- (m_{1,k'+1}x_{k'+1}+\dots + m_{1,k}x_{k})\\
b_2- (m_{2,k'+1}x_{k'+1}+\dots+ + m_{2,k}x_{k})\\
\vdots \\
b_{k'}- (m_{k',k'+1}x_{k'+1}+\dots+ m_{k',k}x_{k})
\end{bmatrix} .$$

Equivalently, set $$ M'_j:=\begin{bmatrix}
m_{1,1} & m_{1,2} &\dots & m_{1,j-1} & b_1- (m_{1,k'+1}x_{k'+1}+\dots + m_{1,k}x_{k}) & \dots & m_{1,k'} \\
m_{2,1} & m_{2,2} &\dots & m_{2,j-1} & b_2- (m_{2,k'+1}x_{k'+1}+\dots+ + m_{2,k}x_{k}) & \dots & m_{2,k'} \\
\vdots & \vdots & \vdots & \vdots & \vdots & & \vdots \\
m_{k',1} & m_{k',2} &\dots & m_{k',j-1}& b_{k'}- (m_{k',k'+1}x_{k'+1}+\dots+ m_{k',k}x_{k})
& \dots & m_{k',k'}
\end{bmatrix}$$
we have that, for any $j\in [1,k']$
\begin{equation}\label{Cramer2}x_j= \frac{\det(M'_j)}{\det(M')}.\end{equation}
Now assume that $R$ is a field and that the system \eqref{sys1} admits at least one solution. Then, given a maximal square nonsingular submatrix $M'$ of $M$, and fixed $x_{k'+1},\dots, x_k$ in $R$, we have that $x_1, x_2, \dots, x_k$ are a solution of the system \eqref{sys1} if and only if $x_1,\dots, x_{k'}$ satisfy the relation \eqref{Cramer2}.
\end{thm}
Now we are ready to prove our Theorem \ref{main}.
\proof
Let $S=\{x_1,\dots,x_k\}$ be a subset of $\mathbb{Z}_m\rtimes_{\varphi} H\setminus\{0\}$ where all the prime factors of $m$ are larger than $k!$.

Given a subset $X$ of $S$ that sum to zero following its ordering $x_{\omega(1)},x_{\omega(2)},\dots,x_{\omega(t-1)},x_{\omega(t)}$, we have the equation, \begin{equation}\label{system} x_{\omega(1)}+ x_{\omega(2)}\dots x_{\omega(t-1)} + x_{\omega(t)}=0.\end{equation}

Here $S$ is a solution $(\pi_{\mathbb{Z}_m}(x_1),\dots,\pi_{\mathbb{Z}_m}(x_k))^T$ to the system of equations derived from \eqref{system} by varing the subset $X$ of $S$ and by considering the projection over $\mathbb{Z}_m$.
In the following, for simplicity, we will denote $\pi_{\mathbb{Z}_m}(x_j)$ by $x_j^1$ and, using this notation, we are considering the system:
\begin{equation}\label{system2} x^1_{\omega(1)}+ \varphi_{\pi_{H}(x_{\omega(1)})}x^1_{\omega(2)}+\dots + \varphi_{\pi_{H}(x_{\omega(1)})\cdots \pi_{H}(x_{\omega(t-1)})}x^1_{\omega(t)}=0\end{equation}
where, using the notation of the previous sections, and since $\varphi$ satisfies property $(*)$, we set $\varphi(a)(x)=\varphi_a x$ for some constant $\varphi_a\in \{1,-1\}$.
Denoted by $M=(m_{i,j})$ the $(k!)\times k$ matrix of the coefficients of this system, we note that, since $\varphi$ satisfies property $(*)$, the coefficients of $M$ belong to $\{-1,0,1\}$. It follows that the determinant of any square submatrix of $M$ (that is big at most $k\times k$) is a sum of at most $k!$ terms that belong to $\pm 1$. Therefore, if all the prime factors of $m$ are larger than $k!$, a square submatrix $M'$ of $M$ is nonsingular in $\mathbb{Z}_{m}$ if and only if it is nonsingular in $\mathbb{Q}$.

Since $\mathbb{Z}_m$ is a commutative ring, once we have chosen a maximal square nonsingular (over $\mathbb{Z}_m$) $k'\times k'$ submatrix $M'$ of $M$, because of Theorem \ref{Cramer}, any solution of the system \eqref{system} satisfies the relations \eqref{Cramer2}
$$x_j^1= \frac{\det(M'_j)}{\det(M')}=\frac{\sum_{i=1}^{k'} (-1)^{j+i+1} \det(M'_{i,j})(m_{i,k'+1}x_{k'+1}^1+\dots + m_{i,k}x_{k}^1)}{\det(M')}$$
for any $j\leq k'$, and where we have assumed, without loss of generality that $M'$ is the $k'\times k'$ matrix obtained by considering the first $k'$ rows and the first $k'$ column of $M$ and where we consider the determinants in $\mathbb{Z}_m$.

Here $M'$ is a maximal square nonsingular $k'\times k'$ submatrix of $M$  also over $\mathbb{Q}$. Hence, following the proof of Theorem 4.12 of \cite{CDOR} we consider the set $S'$ of $\mathbb{Q}\rtimes_{\varphi} H\setminus\{0\}$ whose first components (in $\mathbb{Q}$) are given by
$$z_{j}^1=\begin{cases}
\frac{\sum_{i=1}^{k'} (-1)^{j+i+1} \det((M'_{i,j})_{\mathbb{Q}})(m_{i,k'+1}\widetilde{x_{k'+1}^1}+\dots + m_{i,k}\widetilde{x_{k}^1})}{\det((M')_{\mathbb{Q}})}\mbox{ if } j\leq k';\\
\widetilde{x_j^1} \mbox{ otherwise.}
\end{cases}
$$
where we denote by $\widetilde{x}$ the smallest non-negative integer that is equivalent to $x$ modulo $m$ and whose second ones (in $H$) are the same of $S$. Due to Theorem \ref{Cramer}, the first components of $S'$ satisfy system \eqref{system2}. Also, following again the proof of Theorem 4.12 of \cite{CDOR}, it is easy to see that the elements of $S'$ are distinct and non-zero.
This means that $S$ is Freiman weakly homomorphic to the subset $S'$ of $\mathbb{Q}\rtimes_{\varphi} H\setminus\{0\}$.

We note that by projecting both the numerators and the denominators of the first components of $S'$ in $\mathbb{Z}_m$ we obtain $S$. Hence, for any subset of $S'$ which sums to zero in a given order also the corresponding subset of $S$ sums to zero in the corresponding order. This means that $S$ and $S'$ are Freiman weakly isomorphic.

Finally, by clearing the denominators of the elements of the set $S'$, we prove that $S$ is Freiman weakly isomorphic to a subset of $\mathbb{Z}\rtimes_{\varphi} H\setminus\{0\}$.
\endproof
Here we note that, proceeding as in \cite{NK} and rectify the set $T \cup (-T)\cup \{0\}$, where $T$ is the set of first coordinates of elements of $S$ we could get a shorter proof of Theorem \ref{main} but with a weaker bound (i.e. $(k-1)^{2k+1}$).

If we consider direct products of type $\Z_m \times H$, we can slightly improve the result of Theorem \ref{main}. Indeed in this case the coefficients of $M$ belong to $\{0,1\}$ and hence, with essentially the same proof of Theorem \ref{main}, we can state the following theorem.

\begin{thm}\label{main2}
Let $H$ be a finite group. Then any subset $S$ of size $k$ of $\mathbb{Z}_m\times H\setminus\{0\}$ is Freiman weakly isomorphic to a subset $S'$ of size $k$ of  $\mathbb{Z} \times H\setminus\{0\}$ provided that the prime factors of $m$ are all greater than $k!/2$.
\end{thm}
\section{Sequencing small sets}
\subsection{Semidirect products}
The idea of this subsection is to adapt the rectification procedure of \cite{NK} by using weak Freiman isomorphism to obtain sequenceability in semidirect products (and, in particular in dihedral groups) and to improve the bounds of \cite{NK} in the case of cyclic groups. Indeed Kravitz proved that
\begin{thm}\label{SeqZ}
Let $S$ be a subset of $\mathbb{Z}\setminus\{0\}$. Then $S$ is sequenceable.
\end{thm}
As a consequence, he proved that any subset of size $k$ of $\mathbb{Z}_m\setminus\{0\}$ is sequenceable provided that the prime factors of $m$ are larger than $(k-1)^{k+1}$. The first application of the weak Freiman isomorphism method we present is that we can slightly improve this result as follows.
\begin{thm}\label{impr}
Let $S$ be a set of size $k$ of $\mathbb{Z}_m\setminus\{0\}$. Then $S$ is sequenceable provided that the prime factors of $m$ are larger than $k!/2$.
\end{thm}
\proof
Under this hypothesis, $S=\{x_1,\dots,x_k\}$ is Freiman weakly isomorphic to a subset $S'=\{\phi(x_1),\dots,\phi(x_k)\}$ of $\mathbb{Z}$ due to Theorem \ref{main2}.
Due to Theorem \ref{SeqZ}, $S'$ is sequenceable and hence there exists $\omega$ so that $(\phi(x_{\omega(1)}),\dots,\phi(x_{\omega(k)}))$ is a sequencing.

Here we prove that $(x_{\omega(1)},\dots,x_{\omega(k)})$ is also a sequencing.
Indeed, $s_i=s_j$ means $x_{\omega(i+1)}+x_{\omega(i+2)}+\dots+x_{\omega(j)}=0$, it  would imply that $\phi(x_{\omega(i+1)})+\phi(x_{\omega(i+2)})+\dots+\phi(x_{\omega(j)})=0$. This is not possible because the partial sums $s_i$ and $s_j$ in the sequencing of $S'$ are distinct. It follows that also $(x_{\omega(1)},\dots,x_{\omega(k)})$ is a sequencing.
\endproof
Finally, we apply this procedure to the case of semidirect products. Following \cite{CDOR}, we will call {\em strongly sequenceable} a group whose subsets that do not contain zero are all sequenceable.
First of all we prove the following result for $\Z\rtimes_{\varphi} H$:
\begin{thm}\label{SeqDZ}
Let $H$ be an abelian strongly sequenceable group and let $S$ be a subset of $\Z\rtimes_{\varphi} H\setminus \{0\}$. Then $S$ is sequenceable.
\end{thm}
\proof
We recall that, for any $h\in H$ we have that either $\varphi(h)=id \mbox{ or }\varphi(h)(x)=-x$ and that $\varphi(0)=id$. Here we divide $S$ into the following four sets:

\begin{itemize}
\item $Z:=\{(0,g)\in S: \varphi(g)=id\}$;
\item $P:=\{(x,h)\in S: \varphi(h)=id \mbox{ and }x>0\}$;
\item $N:=\{(x,h)\in S: \varphi(h)=id \mbox{ and }x<0\}$;
\item $F:=\{(y,f)\in S; \varphi(f)=-id\}$.
\end{itemize}
Also, we may assume that the size of $F$ is $\ell\geq1$ otherwise $S$ is morally a subset of $\Z\times H$. More precisely, $S$ is Freiman isomorphic to a subset of $\Z\times H$ and is then sequenceable due to the result of \cite{NK} and because $H$ is abelian.

Now we consider three cases:

CASE $1$: $P=N=Z=\emptyset$.
First, we assume that the projection of $S$ over $\mathbb{Z}$ is constant. In this case, the projection of $S$ over $H$ is a set. Here if $\pi_H(S)$ does not contain $0$, $S$ is sequenceable since $\pi_H(S)$ is sequenceable. On the other hand, if we have an element of type $(y,0)\in S$, we can consider a sequencing
$$((y_2,f_2),(y_3,f_3),\dots,(y_{\ell},f_{\ell}))$$
of $S'=S\setminus \{(y,0)\}$ and place $(y,0)$ at the beginning of this ordering. It is easy to check that
$$((y,0),(y_2,f_2),(y_3,f_3),\dots,(y_{\ell},f_{\ell}))$$
is a sequencing of $S$.

Therefore, we may assume that the projection over $\mathbb{Z}$ of $F$ contains at least two different elements (i.e. not all the $y_i$ are equal).
In this case, we order the elements of $F$ as follows
$((y_2,f_2),(y_{3},f_3),\dots,(y_{\ell-1},f_{\ell-1}),(y_{\ell},f_{\ell}),(y_1,f_1))$
where $$\dots\leq y_7\leq y_5\leq y_3\leq y_1\leq y_2\leq y_4\leq y_6\leq \dots$$
and set
$$F_1:=\{(y_i,f_i)\in F: i>1,\ y_i=y_{i-1}=y_1\}$$
we have that $F_1$ is either empty or, denoted by $\ell_1=|F_1\cup\{(y_1,f_1)\}|$, $(y_2,f_2),\dots, (y_{\ell_1}, f_{\ell_1})$ is a sequencing of $F_1$.

We note that
\begin{equation*}\sum_{x\in F_1} x=\sum_{i=2}^{\ell_1}(y_i,f_i)=(y_2-y_3+y_4-y_5+\dots,\sum_{i=2}^{\ell_1}f_{i})=\end{equation*}
\begin{equation}\label{SumIndependent}(y_1-y_1+y_1-y_1+\dots,\sum_{i=2}^{\ell_1} f_{i})=\left(\frac{y_1+(-1)^{\ell_1} y_1}{2},\sum_{i=2}^{\ell_1} f_{i}\right)\end{equation}
does not depend on the ordering of the elements of $F_1$ (the second component is invariant since $H$ is abelian).
It follows that, up to interchange $(y_1,f_1)$ with any elements in $F_1$, we may assume that:
\begin{equation}\label{SumNonZero}\sum_{x\in F_1} x\not=0.\end{equation}
Here we have that $(y_{\ell},f_{\ell})\not \in F_1$ since otherwise any elements of $F$ would have $y_1$ as a first component. On the other hand, note that $y_{\ell}$ could equal $y_1$.

If $|F_1|\geq 2$, we have at least two choices for $(y_1,f_1)$ so that Equation \eqref{SumNonZero} holds true. At least one of this choice is such that $(y_1,f_1)+(y_{\ell},f_{\ell})\not=0$.

If we have that $|F_1|=0$ then we must have $y_1\not=y_{\ell}$ and hence also in this case we have that $(y_1,f_1)+(y_{\ell},f_{\ell})\not=0$.
Finally, if $|F_1|=1$ (i.e. $F_1=\{(y_2,f_2)\}$), Equation \eqref{SumNonZero} is always satisfied and we can guarantee that $(y_1,f_1)+(y_{\ell},f_{\ell})\not=0$ up to interchanging $(y_1,f_1)$ with $(y_2,f_2)$.

In any case, in the following, we will assume that
$$(y_1,f_1)+(y_{\ell},f_{\ell})\not=0.$$

We recall that $(a_1,a_2,\dots,a_{|S|})$ is a sequencing of $S$ if the partial sums $s_i\not=s_j$ whenever $i\not=j$ and $(i,j)$ are not $(0,|S|)$, which is equivalent to require that \begin{equation}\label{cond2}a_{i+1}+ a_{i+2}+ \dots + a_j\not=0.\end{equation}
Here, since $\varphi$ is a homomorphism from $H$ to $Aut(\mathbb{Z})$, to have zero-sum on the $H$ component, a subset $X$ of consecutive elements must contain an even number of elements of $F$. Indeed, if $X$ contains an odd number of elements of $F$ we would have that
$\varphi(\pi_H(\sum_{x\in X} x))=-id$ that implies $\sum_{x\in X} x\not=0$.

These subsets $X$ are of one of the following forms
\begin{itemize}
\item[a)] $\{(y_{2k},f_{2k}),(y_{2k+1},f_{2k+1}),\dots,(y_{2t},f_{2t}),(y_{2t+1},f_{2t+1})\}$.
Here the $\mathbb{Z}$ component of the sum is
$$(y_{2k}-y_{2k+1})+\dots +(y_{2t}-y_{2t+1})$$
which, if $(y_{2t+1},f_{2t+1})\not\in F_1$ is non-zero since
$y_{2k}-y_{2k+1}\geq 0,\dots, y_{2t}-y_{2t+1}\geq 0$ where at least one inequality is strict.
On the other hand, if $$\{(y_{2k},f_{2k}),(y_{2k+1},f_{2k+1}),\dots,(y_{2t},f_{2t}),(y_{2t+1},f_{2t+1})\}\subseteq F_1,$$ then the $H$ component of the sum is non-zero.
\item[b)] $\{(y_{2k-1},f_{2k-1}),(y_{2k},f_{2k}),\dots,(y_{2t-1},f_{2t-1}),(y_{2t},f_{2t})\}$.
Here the $\mathbb{Z}$ component of the sum is
$$(y_{2k-1}-y_{2k})+\dots +(y_{2t-1}-y_{2t})$$
which, if $(y_{2t},f_{2t})\not\in F_1$ is non-zero since
$y_{2k-1}-y_{2k}\leq 0,\dots, y_{2t-1}-y_{2t}\leq 0$ where at least one inequality is strict.
On the other hand, if $$\{(y_{2k-1},f_{2k-1}),(y_{2k},f_{2k}),\dots,(y_{2t-1},f_{2t-1}),(y_{2t},f_{2t})\}\subseteq F_1,$$ then the $H$ component of the sum is non-zero.

\item[c)] $\{(y_{\ell-2k},f_{\ell-2k}),(y_{\ell-2k+1},f_{\ell-2k+1}),\dots,(y_{\ell-1},f_{\ell-1}),(y_{\ell},f_{\ell}),(y_{1},f_{1})\}$.
Here we may assume that $k\geq 1$ otherwise the sum is non-zero by construction. Then, the $\mathbb{Z}$ component of the sum is
$$(y_{\ell-2k}-y_{\ell-2k+})+\dots +(y_{\ell}-y_{1})$$
which is non-zero since
$y_{\ell-2k}-y_{\ell-2k+1},\dots, y_{\ell}-y_{1}$ are all either positive or negative and at least one inequality is strict.
\end{itemize}
CASE $2$: $P=N=\emptyset$ but $Z\not=\emptyset$. Also in this case we may assume that the projection over $\mathbb{Z}$ of $F$ contains at least two different elements (i.e. not all the $y_i$ are equal) otherwise we can prove, proceeding as in case $1$ that $S$ is sequenceable even if $0\in \pi_{H}(S)$ (because its projection over $H$ is a set).

Then we first order the elements of $F$ as follows
$$((y_1,f_1),(y_2,f_2),(y_{3},f_3),\dots,(y_{\ell-1},f_{\ell-1}),(y_{\ell},f_{\ell}))$$
where $$\dots\leq y_7\leq y_5\leq y_3\leq y_1\leq y_2\leq y_4\leq y_6\leq \dots$$
and where, set
$$F_1:=\{(y_i,f_i)\in F: i\geq 1,\ y_i=y_{i-1}=y_1\}.$$

Note that here $F_1$ contains also $(y_1,f_1)$ and we have not assumed any hypothesis about its ordering, yet. Indeed we will reorder $F_1$ together with $Z$.

We note that, as seen in Equation \eqref{SumIndependent}, the sum $\sum_{x\in F_1\cup Z} x$ does not depend on the ordering. Now, if $\sum_{x\in F_1\cup Z} x=0$ we chose $(0,g_1)\in Z$ and we define $Z_1=Z\setminus \{(0,g_1)\}$ otherwise we set $Z_1=Z$.
In this way, we have that $\sum_{x\in F_1\cup Z_1} x\not=0$ and, since the projection of $F_1\cup Z_1$ over $H$ is a set it admits a sequencing (even if $0\in \pi_H(F_1\cup Z_1)$):
$$C=((z_1,d_1),\dots,(z_{|F_1|+|Z_1|},d_{|F_1|+|Z_1|})).$$
Now, if $Z=Z_1$, we consider the following ordering on $S$
$$(C, (y_{|F_1|+1},f_{|F_1|+1}), (y_{|F_1|+2},f_{|F_1|+2}),\dots  (y_{\ell},f_{\ell}))$$
and, if $Z=Z_1\dot{\cup} \{(0,g_1)\}$, we consider
$$(C, (y_{|F_1|+1},f_{|F_1|+1}), (y_{|F_1|+2},f_{|F_1|+2}),\dots  (y_{\ell},f_{\ell}),(0,g_1)).$$
Also in this case we need to check the sums over subsets of consecutive elements which contain an even number of elements of $F$. The projection of these subsets over $\mathbb{Z}$ has the same form of either case 1.a or case 1.b and is therefore non-zero.

CASE $3$: $P\cup N\not=\emptyset$. In this case, we may assume (up to change the signs of all the $\mathbb{Z}$ components) that $P\not=\emptyset$.

We order the elements of $F$ as follows
$((y_1,f_1),(y_2,f_2),(y_{3},f_3),\dots,(y_{\ell-1},f_{\ell-1}),(y_{\ell},f_{\ell}))$
where $$\dots\leq y_7\leq y_5\leq y_3\leq y_1\leq y_2\leq y_4\leq y_6\leq \dots$$
and where, set
$$F_1:=\{(y_i,f_i)\in F: i>1,\ y_i=y_{i-1}=y_1\}$$
we have that either $F_1$ is empty or $$\sum_{x\in F_1} x\not=0$$
and that $(y_2,f_2),\dots, (y_{|F_1|}, f_{|F_1|})$ is a sequencing of $F_1$.

As in case 1, also here this ordering exists because we can interchange $(y_1,f_1)$ with any elements in $F_1$ and because the sum $\sum_{x\in F_1} x$ does not depend on the ordering of its elements.

Since $H$ is strongly sequenceable, the set $Z$ admits the sequencing $((0,g_1),\dots, (0,g_z))$ where $z=|Z|$. We note that, since $H$ is abelian, also the sum $\sum_{i=1}^z g_i$ does not depend on the ordering of the elements of $Z$. Now, if $\sum_{i=1}^z g_i=0$ we define the set $Z_1$ as $Z\setminus \{(0,g_z)\}$ and the sequence $C=((y_1,f_1),(0,g_z))$ otherwise we set $Z_1=Z$ and $C=((y_1,f_1))$.
In both cases, we define $z'=|Z_1|$.

Also, set $p=|P|$ and $n=|N|$, we fix a generic order of $N$
$$((x_1,h_1),(x_2,h_2),\dots,(x_{n},h_n))$$
and a generic order of $P$
$$((x_{n+1},h_{n+1}),\dots,(x_{n+p},h_{n+p})).$$
Then we consider the following ordering on $S$:
\begin{equation*}((0,g_1),\dots, (0,g_{z'}),(x_1,h_1),(x_2,h_2),\dots,(x_{n},h_n),C\end{equation*}
\begin{equation}\label{orderD}
(x_{n+1},h_{n+1}),\dots,(x_{n+p},h_{n+p}),(y_2,f_2),(y_3,f_3),\dots,(y_\ell,f_{\ell})).
\end{equation}
In the following, when it will be not important to specify if an element of $S$ belongs to $Z,N,P$ or $F$, we will denote the ordering of Equation \eqref{orderD} as:
$$((z_1,d_1),(z_2,d_2),\dots,(z_{|S|},d_{|S|})).$$

Here the subsets $X$ that we need to check are of one of the following forms
\begin{itemize}
\item[a)] $\{(z_k,d_k),(z_{k+1},d_{k+1}),\dots, (z_t,d_t)\}$ which is contained either in $Z_1\cup N$ or in $(Z\setminus Z_1)\cup P$. Since $Z_1$ has distinct partial sums we may assume $X$ is not contained in $Z_1$ and hence, in both of these cases is clear that the sum is non-zero in the $\mathbb{Z}$ component.
\item[b)] $$\{(z_k,d_k),(z_{k+1},d_k),\dots, (x_n,h_n),C,(x_{n+1},h_{n+1}),\dots,$$ $$(x_{n+p},h_{n+p}),(y_2,f_2),(y_3,f_3),\dots, (y_{2t},f_{2t})\}.$$ Here the $\mathbb{Z}$ component of the sum is
$$(z_k+\dots+x_n)+\pi_{\mathbb{Z}}\left(\sum_{x\in C}x\right)-(x_{n+1}+\dots+x_{n+p})-y_2+\dots  +(y_{2t-1}-y_{2t}).$$
This sum is non-zero because $\pi_{\mathbb{Z}}\left(\sum_{x\in C}x\right)=y_1$ and $(z_k+\dots+x_n)\leq 0$, $-(x_{n+1}+\dots+x_{n+p})<0$ and $y_1-y_2\leq 0, y_3-y_4\leq 0,\dots, y_{2t-1}-y_{2t}\leq 0$.

\item[c)] $\{(x_k,h_k),(x_{k+1},h_{k+1}),\dots, (x_{n+p},h_{n+p})$ $,(y_2,f_2),(y_3,f_3),\dots,(y_{2t},f_{2t}),(y_{2t+1},f_{2t+1})\}$ for $k \geq n+1$.
Here the $\mathbb{Z}$ component of the sum is
$$(x_{k}+\dots+x_{n+p})+(y_2-y_3)+\dots +(y_{2t}-y_{2t+1}).$$
This sum is non-zero because $(x_{k}+\dots+x_{n+p})>0$ and $y_2-y_3\geq 0,\dots, y_{2t}-y_{2t+1}\geq 0$.

Note that also the case  $$\{(0,g_z),(x_{n+1},h_{n+1}),\dots, (x_{n+p},h_{n+p}),(y_2,f_2),(y_3,f_3),\dots,(y_{2t},f_{2t}),(y_{2t+1},f_{2t+1})\}$$
can be treated in the same way.
\item[d)] $\{(y_{2k},f_{2k}),(y_{2k+1},f_{2k+1}),\dots,(y_{2t},f_{2t}),(y_{2t+1},f_{2t+1})\}$.
Here we verify that the sum is non-zero proceeding as in case 1.a.

\item[e)] $\{(y_{2k-1},f_{2k-1}),(y_{2k},f_{2k}),\dots,(y_{2t-1},f_{2t-1}),(y_{2t},f_{2t})\}$.
Here we verify that the sum is non-zero proceeding as in case 1.b.
\end{itemize}
\endproof

Applying Theorem \ref{main}, with the same proof of Theorem \ref{impr}, we obtain the following theorem.
\begin{thm}
Let $H$ be a strongly sequenceable finite abelian group, and let $\varphi$ satisfy property $(*)$. Then any subset $S$ of size $k$ of $\mathbb{Z}_m\rtimes_{\varphi} H\setminus\{0\}$ is sequenceable provided that the prime factors of $m$ are all greater than $k!$.
\end{thm}
As a corollary, we find the following result on dihedral groups.
\begin{cor}
Any subset $S$ of size $k$ of $D_{2m}\setminus\{0\}$ is sequenceable provided that the prime factors of $m$ are all greater than $k!$.
\end{cor}
\subsection{Dicyclic groups}
In this subsection, we show that any subset of size $k$ of the dicyclic group $\Dic_{m}$ is sequenceable, provided that the prime factors of $m$ are large enough.
\begin{defn}\label{Dic}
We denote by $\Dic_m$ the group of order $4m$ whose generators are the elements $s$ and $r$ which satisfy the following relations:
\begin{itemize}
\item[1)] $r+s=s-r$;
\item[2)] $4s=0$;
\item[3)] $2s=mr$.
\end{itemize}
\end{defn}
\begin{rem}\label{SemidirectZZ4}
We observe that the group whose generators are the elements $s$ and $r$ which satisfy relations $1$ and $2$ of Definition \ref{Dic}, is the semidirect product $\mathbb{Z}\rtimes_{\varphi} \mathbb{Z}_4$ where $\varphi(1)=-id$. We note that, due to Theorem \ref{SeqDZ}, $\mathbb{Z}\rtimes_{\varphi} \mathbb{Z}_4$ is strongly sequenceable.
\end{rem}
Our goal is now to show that subsets of dicyclic groups are weak Freiman isomorphic to subsets of this semidirect product. We first need the following technical lemma.
\begin{lem}\label{tech}
Let $S$ be a subset of $\mathbb{Z}_{2m}$ where the prime factors of $m$ are larger than $k^k$. Then there exists an automorphism $\alpha: \mathbb{Z}_{2m} \rightarrow \mathbb{Z}_{2m}$ such that 
$$\alpha(S)\subseteq \left\{-\lfloor\frac{m}{k}\rfloor,-\lfloor\frac{m}{k}\rfloor+1,\dots, \lfloor\frac{m}{k}\rfloor\right\} \cup \left(m+ \left\{-\lfloor\frac{m}{k}\rfloor,-\lfloor\frac{m}{k}\rfloor+1,\dots, \lfloor\frac{m}{k}\rfloor\right\}\right).$$
\end{lem}
\proof
We observe that $m$ must be odd and hence $\mathbb{Z}_{2m}\cong \mathbb{Z}_2 \times \mathbb{Z}_m$. The isomorphism between these groups is given by the Chinese remainder theorem which maps $x$ to $(x \pmod{2}, x\pmod{m})$ and its inverse maps $(x_1,x_2)$ into $mx_1+(m+1)x_2$.
It follows that, in order to prove the thesis of the theorem, it suffices to show that given a subset of $\mathbb{Z}_m$, say $S'=\{x_1,x_2,\dots,x_k\}$, $\mathbb{Z}_m$ admits an automorphism $\alpha_1: \mathbb{Z}_{m} \rightarrow \mathbb{Z}_{m}$ such that
$$\alpha_1(S')\subseteq \left\{-\lfloor\frac{m}{k}\rfloor,-\lfloor\frac{m}{k}\rfloor+1,\dots, \lfloor\frac{m}{k}\rfloor\right\}.$$

Let $p$ be the smallest prime factor of $m$. Then applying Minkowski's theorem as done in the proof of Theorem 3.1\footnote{ We can apply that same argument since the matrix of Minkowski's theorem is lower triangular and its determinant is $\pm 1$ independently on the coefficients $x_i/m$.} of \cite{BLR}, and since $p>k^k$, we find that there exist integers $a,a_1,\dots, a_k$ such that
\begin{equation}\label{Mink1}1\leq |a|<p,\end{equation}
\begin{equation}\label{Mink2}\left|\frac{a x_i}{m} - a_i\right| \leq \frac{1}{k} \qquad \text{for} \qquad 1 \leq i \leq k.\end{equation}
It follows, from Equation \eqref{Mink2}, that there exist $$y_1,y_2,\dots, y_k \in \left\{-\lfloor m/k \rfloor, -\lfloor m/k \rfloor +1, \dots, \lfloor m/k \rfloor \right\}$$ such that $y_i \equiv a x_i \pmod{m}$.

The thesis follows since, due to Equation \eqref{Mink1}, $a\not \equiv 0 \pmod {q}$ for any prime factor $q$ of $m$ and hence the map $\alpha_1: x\rightarrow ax$ is the required automorphism of $\mathbb{Z}_m$.
\endproof
\begin{prop}\label{IsoDic}
Any subset of size $k$ of the dicyclic group $\Dic_m$ is $(k-1)$-weak Freiman isomorphic to a subset of the semidirect product $\mathbb{Z}\rtimes_{\varphi} \mathbb{Z}_4$ (where $\varphi(1)=-id$) provided that the prime factors of $m$ are larger than $k^k$.
\end{prop}
\proof
Due Definition \ref{Dic}, any element of $\Dic_m$ can be written, uniquely, as $\lambda_1s+\lambda_2r$ where $\lambda_1\in [0,3]$ and $\lambda_2\in [-\lfloor m/2\rfloor,\lfloor m/2\rfloor]$. Let now consider a set $S=\{x_1,\dots,x_k\}$ of size $k$ and let us write (as indicated above) $x_i=\lambda_{i,1}s+\lambda_{i,2}r$.

Since the group generated by $r$ is isomorphic to $\mathbb{Z}_{2m}$, due to Lemma \ref{tech}, we can find an isomorphism $\sigma$ of $\Dic_m$ such that
\begin{itemize}
\item $\sigma(s)=s$;
\item $\sigma(\lambda_{i,2}r)=\mu_{i,1}s+\mu_{i,2}r$ where $\mu_{i,1}\in \{0,2\}$ and $$\mu_{i,2}\in \left\{-\lfloor m/k \rfloor, -\lfloor m/k \rfloor +1, \dots, \lfloor m/k \rfloor \right\}.$$
\end{itemize}
This means that the elements $\{y_1,\dots,y_k\}$ of $\sigma(S)$ can be written as
$y_i=\delta_{i,1}s+\delta_{i,2}r$ where $\delta_{i,1}\in [0,3]$ and $\delta_{i,2}\in \left\{-\lfloor m/k \rfloor, -\lfloor m/k \rfloor +1, \dots, \lfloor m/k \rfloor \right\}.$

Now, to show that $S$ is $(k-1)$-weak Freiman isomorphic to a subset of the semidirect product $\mathbb{Z}\rtimes_{\varphi} \mathbb{Z}_4$, it suffices to prove this statement for $\sigma(S)$.
For this purpose we consider the set $S'=\{z_1,z_2,\dots, z_k\}$ where $z_i=(\delta_{i,2},\delta_{i,1})$ and let us define $\phi: \sigma(S)\rightarrow S'$ as $\phi(y_i)=z_i$. Let us consider an ordered subset $Y$ of size $\ell\leq k-1$ of $S$ which sums to zero. Without loss of generality we may assume $Y=(y_1,\dots,y_{\ell}).$
Since its sum is zero, we have that
$$0=y_1+\dots+y_{\ell}= \sum_{i=1}^{\ell}\left (\delta_{i,1}s+\delta_{i,2}r\right)=$$
\begin{equation}\label{somma} \sum_{i=1}^{\ell}\left (\delta_{i,1}\right)s+ \sum_{i=1}^{\ell}\left ((-1)^{\sum_{j=i+1}^\ell \delta_{j,1}}\delta_{i,2}\right)r.\end{equation}
Since $\delta_{i,2}\in \left\{-\lfloor m/k \rfloor, -\lfloor m/k \rfloor +1, \dots, \lfloor m/k \rfloor \right\},$ we have that
$$\left|\sum_{i=1}^{\ell}\left ((-1)^{\sum_{j=i+1}^\ell \delta_{j,1}}\delta_{i,2}\right)\right|\leq \sum_{i=1}^{\ell}|\delta_{i,2}|\leq \frac{m\ell}{k}<m.$$
It follows that Equation \eqref{somma} is satisfied if and only if $$\sum_{i=1}^{\ell}\left ((-1)^{\sum_{j=i+1}^\ell \delta_{j,1}}\delta_{i,2}\right)=0 \qquad \text{and} \qquad \sum_{i=1}^{\ell}\delta_{i,1}\equiv 0\pmod{4}.$$
Since also in $\mathbb{Z}\rtimes_{\varphi}\mathbb{Z}_4$ the sum respects relation $1$ of Definition \ref{Dic}, we also have that 
$$ z_1+\dots+z_{\ell}=\left(\sum_{i=1}^{\ell}\left ((-1)^{\sum_{j=i+1}^\ell \delta_{j,1}}\delta_{i,2}\right),\sum_{i=1}^{\ell}\delta_{i,1}\right).$$
Therefore, $0=y_1+\dots+y_{\ell}$ implies $0=z_1+\dots+z_{\ell}$.

On the other hand, since $\Dic_m$ is defined extending the relations of  $\mathbb{Z}\rtimes_{\varphi}\mathbb{Z}_4$, it is clear that also the converse holds true and hence $\phi$ is a $(k-1)$-weak isomorphism between $\sigma(S)$ and $S'$. As observed above, this implies the thesis.
\endproof
As observed in Remark \ref{SemidirectZZ4}, due to Theorem \ref{SeqDZ}, $\mathbb{Z}\rtimes_{\varphi}\mathbb{Z}_4$ is strongly sequenceable. Then,  as a consequence of Proposition \ref{IsoDic} we obtain the following sequenceability result on dicyclic groups.
\begin{thm}
Any subset $S$ of size $k$ of $\Dic_{m}\setminus\{0\}$ is sequenceable provided that the prime factors of $m$ are all greater than $k^k$.
\end{thm}
\section*{Acknowledgements} The authors would like to thank Noah Kravitz for our useful discussions on this
topic and for pointing out that Theorem \ref{SeqDZ} can be generalized to semidirect products and suggesting that our approach can be applied to dicyclic groups. 

The authors were partially supported by INdAM--GNSAGA.

\end{document}